\documentclass[12pt]{article}
\usepackage{longtable}

\usepackage{amssymb}
\usepackage{amsmath}
\usepackage{algorithmic}
\usepackage{graphicx}
\usepackage{amsthm}
\newcommand{\Gr}{Gr\"obner }

\newcommand{\K}{\mathbb{K}}
\newcommand{\N}{\mathbb{N}}
\newcommand{\M}{\mathbb{M}}
\newcommand{\R}{\mathbb{R}}
\newcommand{\X}{\mathbb{X}}

\newcommand{\lm}{\mathop{\mathrm{lm}}\nolimits}
\newcommand{\lt}{\mathop{\mathrm{lt}}\nolimits}

\newcommand{\Id}{\mathop{\mathrm{Id}}\nolimits}

\newcommand{\mindeg}{\mathop{\mathrm{mindeg}}\nolimits}

\def \bg #1 {\begin{tabular}{{#1}}}
\def \nd {\end{tabular}}

\newenvironment{algorithm}[1]{
  \begin{center}
    {\bf Algorithm: #1}\\*
     \begin{tabular}{|p{130mm}|} \hline
} {
 \\ \hline
 \end{tabular}
 \end{center}
}
\newtheorem{definition}{Definition}[section]

\begin{document}

\title{\bf On Computing Janet Bases
for Degree Compatible Orderings}

\author{Vladimir P. Gerdt \\
       Laboratory of Information Technologies\\
       Joint Institute for Nuclear Research\\
       141980 Dubna, Russia \\
       gerdt@jinr.ru
\and
       Yuri A. Blinkov \\
       Department of Mathematics and Mechanics \\
       Saratov State University \\
       410071 Saratov, Russia \\
       blinkovua@info.sgu.ru}
\date{}
\maketitle

\begin{abstract}
We consider three modifications of our involutive algorithm for computing
Janet bases. These modifications are related to degree compatible monomial
orders and specify selection strategies for non-multiplicative
prolongations. By using the standard data base of polynomial benchmarks for
\Gr bases software we compare the modifications and confront them with Magma
that implements Faug\`{e}re's $F_4$ algorithm.
\end{abstract}

\section{Introduction}

In~\cite{GB'98} we designed an algorithmic approach to computing \Gr bases based on
the new notion of involutive monomial division. This notion generalizes the constructive ideas
of French mathematician M.Janet~\cite{Janet} on partition of independent variables
for partial differential equations (PDEs) into multiplicative and non-multiplicative. He applied
this partition to complete systems of PDEs to involutivion. Janet's receipt of the separation of
variables generates a particular involutive division which we called in~\cite{GB'98} Janet
division.

Our algorithms~\cite{GBY'01,G'05}, which improve those in~\cite{GB'98}, shown rather good computational
efficiency when specialized for Janet division. Moreover, for the present, none of other known
involutive divisions demonstrated its computational superiority over Janet division. Though the improved
algorithms allow a user to output a reduced \Gr basis without extra reductions in the final involutive basis, it
is computationally worthwhile to minimize \Gr redundancy in the intermediate and in final
involutive basis. Recently~\cite{GB'05}, we succeeded in such an improvement of Janet division,
and called the improved division as Janet-like. The last division, since it does not partition the variables
into multiplicative and non-multiplicative, is not involutive. But even so, it is very close to
Janet division and
preserves all its computational merits. Quite often, e.g.  for most of benchmarks in the standard
data base~\cite{BM'96,JV}, there is no substantial computational superiority of Janet-like division over
Janet divisions. But for other classes of polynomial problems such as toric ideals this superiority can be
enormous~\cite{GB'05}.

Apart from improvement of the division, there is another important source of optimization in the
involutive algorithms: selection of non-multiplicative prolongations. The last constructions, i.e.,
the products
of intermediate polynomials and their non-multiplicative variables, in the case of their involutive
head reducibility, play in the involutive
algorithms the same role
as $S-$polynomials play in Buchberger's algorithm~\cite{Buch'85}. In the latter case any $S-$polynomial
can be selected at any
step of the algorithm. This enormous arbitrariness, as well as that in the reduction strategy, is
one of the main obstacles
on the way of optimization of Buchberger's algorithm. Only after many years of research some
heuristically good selection
strategies for $S-$polynomials, such as "sugar"~\cite{sugar}, were found.

As to the involutive approach, the reduction sequence is uniquely defined and an admissible
choice of a non-multiplicative prolongation
is subject to certain restrictions~\cite{G'05}. Nevertheless, for examples large enough, one can
choose from many possible prolongations.
For example, in the 7th order cyclic root example~\cite{BM'96,JV} at the intermediate algorithmic
steps there arise several hundreds
prolongations such that any of them can be chosen. By this reason it is important to
investigate a heuristical efficiency of different selection strategies.

In the given paper we present three different selection strategies which as we found
are computationally good. In so going, we restrict
ourselves with degree compatible monomial orders. In practice, this is a reasonable restriction.
It is well-known that
heuristically best way of computing a \Gr basis for an arbitrary order is to compute
it, first,
for a degree compatible order
and then to convert the basis into the desirable one by the FGLM algorithm~\cite{FGLM} or by the \Gr
walk~\cite{GroebnerWalk}.
After some preliminary definitions and conventions (Sect.2) we modify in Sect.3 the involutive
algorithm~\cite{G'05}, as
specialized for Janet division, in accordance to those selection strategies. In Sect.4 we give
experimental comparison
of the three modifications of involutive algorithm on the benchmarks from~\cite{BM'96,JV}. Here we
also show the corresponding timings
for the implementation in Magma~\cite{Magma} of Faug\'ere's $F_4$ algorithm~\cite{F4}.
We conclude in Sect.5.

\section{Preliminaries}
Throughout this paper we shall use the following notations and conventions:
\begin{description}

\item[] $\N_{\geq 0}$ is the set of nonnegative integers. \\[-0.8cm]

\item[] $\X=\{x_1,\ldots,x_n\}$ is the set of polynomial variables.\\[-0.8cm]

\item[] $\R=\K[\X]$ is a polynomial ring over a zero characteristic field $\K$.\\[-0.8cm]

\item[] $\Id(F)$ is the ideal in $\R$ generated by $F\subset \R$.\\[-0.8cm]

\item[] $\M=\{x_1^{i_1}\cdots x_n^{i_n} \mid i_k\in \N_{\geq 0},\
1\leq k\leq n\}$ is the monoid of monomials in $\R$.\\[-0.8cm]

\item[] $\deg_i(u)$ is the degree of $x_i$ in $u\in \M$.\\[-0.8cm]

\item[] $\deg(u)=\sum_{i=1}^n \deg_i(u)$ is the total degree of $u$.\\[-0.8cm]

\item[] $\mindeg(U)=\min\{\, \deg(u_1),\ldots,\deg(u_k)\, \}$ for $U=\{u_1,\ldots,u_k\}\subset \M$.
        \\[-0.8cm]

%\item[] $\maxdeg(U)=\max\{\, \deg(u_1),\ldots,\deg(u_k)\, \}$ for $U=\{u_1,\ldots,u_k\}\subset \M$.
%        \\[-0.8cm]

\item[] $\succ$ is an admissible monomial order such that
$x_1\succ x_2\succ\cdots\succ x_n$ and for $u,v\in \M$ \\ $\deg(u)>\deg(v)\Longrightarrow u\succ v$.\\[-0.8cm]

\item[] If monomial $u$ divides monomial $v$ and
$\deg(u)<\deg(v)$, i.e. $u$ is a proper
divisor of $v$, we shall write $u\sqsubset v$.\\[-0.8cm]

\item[] $\lm(f)$ and $\lt(f)$ are the leading monomial and the
leading term of $f\in
\R\setminus \{0\}$.\\[-0.8cm]

\item[] $\lm(F)$ is the leading monomial set for $F\subset \R\setminus \{0\}$.\\[-0.8cm]

\item[] ${J}$ is Janet division.\\[-0.8cm]

\item[] $M_{J}(u,U)$ is the set of ${J}-$multiplicative
variables for
monomial $u\in U\subset \M$.\\[-0.8cm]

\item[] $NM_{J}(u,U)$ is the set of ${J}-$non-multiplicative
variables of
monomial $u\in U\subset \M$.\\[-0.8cm]

\item[] $J(u,U)$ is the submonoid in $\M$ generated by the power products
of variables in $M_{J}(u,U)$.\\[-0.8cm]

\item[] $u\in U$ ($U\subset \M$ is finite) is a Janet divisor of
$v\in \M$ if $v=u\cdot w$, $w\in J(u,U)$.\\[-0.8cm]

\item[] $HNF_{J}(f,F)$ is the ${J(anet)}-$head normal form of $f\in
\R$ modulo
$F\subset \R$.\\[-0.8cm]

\item[] $NF_{J}(f,F)$ is the ${J(anet)}-$(full) normal form of $f\in
\R$ modulo
$F\subset \R$.

\end{description}

\noindent
The following definitions are taken from~\cite{GB'98}.

\begin{definition}{\em {\em Janet division}. Let $U\subset \M$
be a finite set. For each $0\leq i\leq n$ partition $U$ into groups
labeled by non-negative integers $d_0,\ldots,d_i$~($U=[0]$).
$$
[d_0,d_1,\ldots,d_i]:=\{u\in U \mid d_0=0,d_1=\deg_1(u),\cdots, d_i=\deg_i(u) \}. \label{groups}
$$
Indeterminate $x_i$ is {\em $J$(anet)-multiplicative} for $u\in U$ if
$u\in [d_0,\ldots,d_{i-1}]$ and
$\deg_i(u)=\max\{\deg_i(v) \mid v\in [d_0,\ldots,d_{i-1}]\}$.
}
\end{definition}

\begin{definition}{\em({\em ${J}-$reduction}). Given a monomial order $\succ$, a finite set
$F\in \R\setminus \{0\}$ of polynomials and a polynomial $p\in \R\setminus \{0\}$,
we shall say that:
 \begin{enumerate}
 \renewcommand{\theenumi}{(\roman{enumi})}
 \item $p\in \R$ is {\em ${J}-$reducible} {\em modulo} $f\in F$ if
  $p$ has a term $t=a\,u$ ($a\in \K, u\in \M, a\neq 0$) such that $\lm(f)$ is $J-$divisor of $u$, that is,
  $u=\lm(f)\cdot v$ where $v\in {J}(\lm(f),\lm(F))$.
  It yields the {\em ${J}-$reduction} $p\rightarrow g:=p-(a/lc(f))\,f\cdot v$.\\[-0.75cm]
 \item $p$ is {\em ${J}-$reducible modulo} $F$ if there is $f\in F$ such
  that $p$ is ${J}-$reducible modulo $f$.\\[-0.75cm]
 \item $p$ is {\em in the ${J}-$head normal form modulo $F$} ($p=HNF_{{J}}(p,F)$)
  if $\lm(p)$ has no ${J}-$divisors in $\lm(F)$.\\[-0.75cm]
 \item $p$ is {\em in the ${J}-$normal form modulo $F$} ($p=NF_{{J}}(p,F)$)
  if $p$ is not ${J}-$reducible modulo $F$.
\end{enumerate}
} \label{J-red}
\end{definition}

\begin{definition}{\em({\em ${J}-$basis})
A polynomial set $F$ is called
{\em $J-$autoreduced} if each term in every $f\in F$ has no $J-$divisors
in $\lm(F)\setminus \lm(f)$. A $J-$autoreduced set $F$ is called a {\em Janet basis} of $\Id(F)$ if
\begin{equation}
 (\forall f\in F)\ (\forall x\in NM_J(f,F))\ \ [\ NF_J(f\cdot x,F)=0\ ]\,. \label{J_basis}
\end{equation}
A Janet basis $G$ is called {\em minimal} if for
any other Janet basis $F$ of the same ideal the inclusion
$lm(G)\subseteq \lm(F)$ holds.
} \label{JB}
\end{definition}

A Janet basis is a \Gr one~\cite{GB'98}, though generally not
reduced. However, similarly to a reduced \Gr basis, a monic minimal Janet
basis is uniquely defined by an ideal and a monomial order.
In that follows we deal with minimal Janet bases only and omit the word "minimal".

\section{Modified Involutive Algorithm}

Our Janet division algorithm is given by (cf.~\cite{GBY'01,G'05}):
\vskip 0.3cm
\noindent
\begin{algorithm}{JanetBasis ($F,\prec $)\label{JanetBasis}}
\begin{algorithmic}[1]
\INPUT $F\in \R\setminus \{0\}$, a finite set;\ $\prec$, a degree compatible order
\OUTPUT $G$, a Janet basis of $\Id(F)$
\STATE {\bf choose} $f\in F$ of the minimal $\deg(\lm(f))$
\STATE $G:=\{f\}$
\STATE $Q:=F\setminus G$
\DOWHILE
  \STATE $h:=0$
  \WHILE{$Q\neq \emptyset$\ and $h=0$}
    \STATE {\bf choose} $p\in Q$ with minimal $\lm(p)$ w.r.t. $\succ$
    \STATE $Q:=Q\setminus \{p\}$
    \STATE $h:=NF_J(p,G)$
  \ENDWHILE
  \IF{$h\neq 0$}
    \FORALL{$\{ g\in G \mid \lm(g)\sqsupset\, \lm(h)\}$}
      \STATE $Q:=Q\cup \{g\}$; \ $G:=G\setminus \{g\}$
    \ENDFOR
    \STATE $G:=G\cup \{ h \}$
    \STATE $Q:=Q\cup \{\,g\cdot x \mid g\in G,\,
                                    x\in NM_{J}(\lm(g),\lm(G))\,\}$
 \ENDIF
\ENDDO{$Q \neq \emptyset$}
\RETURN $G$
\end{algorithmic}
\end{algorithm}

\noindent
In its improved form~\cite{G'05}, at the initialization step, i.e., before starting the main loop 4-22
and after its modification in the loop, the set $Q$ is $J-$head reduced modulo $G$. Now, as the
first modification of the above algorithm {\bf JanetBasis} we shall use only partial head reduction of
elements in $Q$ as shown in the following algorithm:

\vskip 0.3cm
\noindent
\begin{algorithm}{JanetBasis I ($F,\prec $)\label{JanetBasis I}}
\begin{algorithmic}[1]
\INPUT $F\in \R\setminus \{0\}$, a finite set;\ $\prec$, a degree compatible order
\OUTPUT $G$, a Janet basis of $\Id(F)$
\STATE {\bf choose} $f\in F$ of the minimal $\deg(\lm(f))$
\STATE $G:=\{f\}$
\STATE $Q:=F\setminus G$
\DOWHILE
  \STATE $S:=\{\,q\in Q \mid \deg(\lm(q))=\mindeg(\lm(Q))\,\}$
  \STATE $P:=\emptyset$;\ \ \ $Q:=Q\setminus S$
  \FORALL{$s\in S$}
    \STATE $S:=S\setminus \{s\}$;\ \ $p:=HNF_J(s,G)$
    \IF{$p\neq 0$}
     \STATE $P:=P\cup \{p\}$
    \ENDIF
  \ENDFOR
  \WHILE{$P\neq \emptyset$}
    \STATE {\bf choose} $p\in P$ with minimal $\lm(p)$ w.r.t. $\succ$
    \STATE $P:=P\setminus \{p\}$;\ \ \ $h:=NF_J(p,G)$
    \FORALL{$\{ g\in G \mid \lm(g)\sqsupset\, \lm(h)\}$}
      \STATE $Q:=Q\cup \{g\}$; \ \ $G:=G\setminus \{g\}$
    \ENDFOR
    \STATE $G:=G\cup \{ h \}$
    \STATE $Q:=Q\cup \{\,g\cdot x \mid g\in G,\,
                                    x\in NM_{J}(\lm(g),\lm(G))\,\}$
  \ENDWHILE
\ENDDO{$Q \neq \emptyset$}
\RETURN $G$
\end{algorithmic}
\end{algorithm}

\vskip 0.3cm
\noindent
In this algorithm at step 5 all the polynomials in $Q$ of the minimal head degree are collected in set $S$
and then are
$J-$head reduced modulo $G$ in the {\bf for} loop 7-12 with the collection of nonzero head reduced polynomials
in $P$ at step 10. Then in the {\bf while} loop 13-21 the polynomial in $P$ with the least leading term
is sequentially selected at step 14 of the loop and inserted in set $G$ after its tail $J-$reduction of
step 15. As well as in algorithm {\bf JanetBasis} the displacement of some polynomials from $G$ to $Q$
done at step 17
provides minimality of the output Janet basis~\cite{GB'98}.

Subalgorithms $HNF_J$ and $NF_J$ which are called at steps 8 and 15 compute the Janet head and the full
normal forms, respectively, in accordance with Definition~\ref{J-red}.

\noindent
Apparently, the modifications done in algorithm {\bf JanetBasis I} in comparison with algorithm
{\bf JanetBasis} do not violate its correctness. As well as for the latter algorithm, when the main
{\bf do-while} loop terminates, the polynomial set $G$ satisfies the conditions~(\ref{J_basis}) in
Definition~\ref{JB}.

We consider now another modification of algorithm {\bf JanetBasis} shown in the form of algorithm
{\bf JanetBasis II}.

\vskip 0.3cm
\noindent
\begin{algorithm}{JanetBasis II ($F,\prec $)\label{JanetBasis II}}
\begin{algorithmic}[1]
\INPUT $F\in \R\setminus \{0\}$, a finite set;\ $\prec$, a degree compatible order
\OUTPUT $G$, a Janet basis of $\Id(F)$
\STATE {\bf choose} $f\in F$ of the minimal $\deg(\lm(f))$
\STATE $G:=\{f\}$
\STATE $Q:=F\setminus G$
\DOWHILE
  \STATE $S:=\{\,q\in Q \mid \deg(\lm(q))=\mindeg(\lm(Q))\,\}$
  \STATE $P:=\emptyset$;\ \ \ $Q:=Q\setminus S$
  \FORALL{$s\in S$}
    \STATE $S:=S\setminus \{s\}$;\ \ $p:=NF_J(s,G)$
    \IF{$p\neq 0$}
     \STATE $P:=P\cup \{p\}$
    \ENDIF
  \ENDFOR
  \STATE $P:={\bf Update}(P,\prec)$
  \FORALL{$p\in P$}
    \FORALL{$\{ g\in G \mid \lm(g)\sqsupset\, \lm(p)\}$}
      \STATE $Q:=Q\cup \{g\}$; \ \ $G:=G\setminus \{g\}$
    \ENDFOR
    \STATE $G:=G\cup \{ p \}$
    \STATE $Q:=Q\cup \{\,g\cdot x \mid g\in G,\,
                                    x\in NM_{J}(\lm(g),\lm(G))\,\}$
  \ENDFOR
\ENDDO{$Q \neq \emptyset$}
\RETURN $G$
\end{algorithmic}
\end{algorithm}

\vskip 0.3cm
\noindent
As distinct from algorithm {\bf JanetBasis I}, the full $J-$normal form is computed at step 8.
Besides,
at step 13 the polynomial set $P$ whose elements are inserted in $G$ in the {\bf for} loop 14-20
is updated in accordance to the below subalgorithm {\bf Update}.

At steps 1 and 4 of the subalgorithm we indicate
two different options for the choice of element $p$: with the highest or with
the lowest leading monomial with respect to the order $\succ$.
In our numerical experiments presented in the next section we apply these two different upgrade
strategies when only one of the indicated choices (highest or lowest) is used in the whole
run of the algorithm.

\vskip 0.3cm
\begin{algorithm}{Update$(P,\succ)$}
\begin{algorithmic}[1]
\INPUT $P\subset \R\setminus \{0\}$, a finite set; $\succ$, an order
\OUTPUT $H\subset \R\setminus \{0\}$, an updated input set
   \STATE {\bf choose} $f\in P$ with the highest/lowest $\lm(f)$ w.r.t. $\succ$
   \STATE $H:=\{f\}$;\ \ \ $P:=P\setminus \{f\}$
   \WHILE{$P\neq \emptyset$}
      \STATE {\bf choose} $p\in P$ with the highest/lowest $\lm(p)$ w.r.t. $\succ$
      \STATE $P:=P\setminus \{p\}$
      \STATE $h:=NF_J(p,H)$
        \IF{$h\neq 0$}
           \STATE $H:=H\cup \{h\}$
        \ENDIF
   \ENDWHILE
   \RETURN $H$
\end{algorithmic}
\end{algorithm}

\noindent
In subalgorithm {\bf Update} an element $f$ in the input polynomial set $P$ (which is
$J-$reduced modulo polynomial set $G$ when the subalgorithm is invoked in the main algorithm)
is chosen at the initialization step 1 with the highest or lowest leading term, depends
on the selection strategy used. After that, in the first run of the {\bf while} loop 3-10 the
other polynomial $p$ in $P$, if any, with the same leading monomial as that in $f$ is $J-$reduced
modulo $f$. In the case of nonzero reduction (when monic $p$ is different from monic $f$) the
normal form obtained is added to $f$ at step 8 to be involved in the further reductions. Then, the
processes of the selection and $J-$reduction of elements in $P$ is repeated until $P$ becomes
empty.

The above described modifications related to certain selection strategies for non-multiplicative
prolongations are easily adapted to the improved version of involutive algorithm~\cite{G'05}.
The improved version avoids useless repeated prolongations and applies the involutive analogues
of Buchberger's criteria for detection of some superfluous reductions. Furthermore, it is
straightforward to include the modifications into the Janet-like division
algorithms~\cite{GB'05}.

\section{Computer Experiments}

The improved version of algorithm {\bf JanetBasis I} was implemented in C on as a part of package
JB~\cite{GBY'01} whose version is also included in the library of Singular~\cite{Singular}
and in C++ as a part of the open source software~GINV~\cite{Ginv}. The last software
implements also algorithm~{\bf JanetBasis II} in its improved version and
for both options in subalgorithm {\bf Update}. For all that Ginv implements also Janet-like
division~\cite{GB'05}.

We examined the three selection strategies of Sect.3 by the standard data base of polynomial
benchmarks~\cite{BM'96,JV} and for degree-reverse-lexicographical monomial order. Some of
the benchmarks are listed in the below table together with the
timings they took for computing \Gr bases. For comparison, we also included the timings
of the last two versions of Magma~\cite{Magma} with the fastest \Gr bases module among all
computer algebra systems.  This is owing to implementation of the Faug\`{e}re $F_4$ algorithm
which rests upon linear algebra for doing multiple reductions in contract to
Buchnerger's algorithm~\cite{Buch'85} or our involutive algorithms which are relayed on
the chains of elementary reductions.

As we noticed in Introduction (Sect.1), the involutive algorithm in its improved
form~\cite{G'05} can output a reduced \Gr basis as an internally fixed part of Janet basis
without any extra reductions. By this reason all software included in the table output
the same bases.

It should be noted that GINV for can also use Janet-like division instead of
Janet division whereas the package JB implements Janet division only. Generally,
intermediate polynomial
sets for Janet-like division
are more compact then those for Janet division. There are whole
classes of interesting multivariate polynomial problems, for example, binomial
toric ideals~\cite{BLR'99} closely related to integer programming~\cite{CoTr'91}
for which Janet-like division leads to enormous gain in comparison with Janet division.
However, as we already mentioned (Sect.1), the difference of two divisions does not manifest
itself significantly for benchmarks in the table.

The timings in the table were obtained on the following machines:
\begin{description}
\item[\ \ \ \ JB:] 2xOpteron-242 (1.6 Ghz) with 4Gb of RAM running under Gentoo Linux 2004.3
with gcc-3.4.2 compiler.
\item[\ \ \ \ GINV:] Turion-3400 (1.8 Ghz) with 2Gb of RAM running under Gentoo Linux 2005.1
with gcc-3.4.4 compiler.
\item[\ \ \ \ Magma:] dual processor Pentium III (1 Ghz) with 2 GB of RAM running under SuSE
Linux 8.0 (kernel 2.4.18-64GB-SMP) with gcc-2.95.3 compiler.
\end{description}
All timings in the table are given in seconds, and (*) shows that the example was not computed because of the
memory overflow.

The 2nd and 3th columns in the table show the results for algorithm {\bf JanetBasis I} whereas
the 4th and 5th
columns shows those for algorithm {\bf JanetBasis II} with the choice of the highest and lowest
option in subalgorithm {\bf Update}, respectively.

One can see rather high stability of the involutive algorithm with respect to three variations
used for the selection strategy. In addition, in all three cases the number of
redistributions between $G$ and $Q$ was experimentally tested to be minimal.
In so doing, we observed that the strategy in algorithm {\bf JanetBasis II}
with the lowest element choice in subalgorithm {\bf Update} leads to a slightly more smooth
growth of the intermediate memory needed than its counterpart with the highest element choice.
The hcyclic8 example in the table explicitly demonstrates this observation.

As to
comparison with Magma, it clearly signals on superiority of the linear algebra based $F_4$
algorithm over our reduction strategy that uses the elementary reduction chains. In our future
research plans there is also improvement of the involutive algorithm by doing reductions by
means of linear algebra.

\newpage
\begin{center}
{\bf Benchmarking}
\begin{scriptsize}
\begin{longtable}{|l|r|r|r|r|r|r|}
\hline
Example & Algorithm I    & Algorithm I      & Algorithm II & Algorithm II & Magma  & Magma \\
        &  (JB)\ \ \ \ \ & (GINV)\ \ \ \   &  high\ (GINV)\ \ & low\ (GINV)\ \    & V2.11-8 & V2.12-17\\
\hline
\endhead
\hline
\endfoot

assur44 & \hfill 10.35 & \hfill 14.20 & \hfill 6.33 & \hfill 6.4 &  \hfill 4.56 & \hfill 4.99\\

butcher8 & \hfill 1.06 & \hfill 1.02 & \hfill 0.38 & \hfill 0.39 &  \hfill 4.68 & \hfill 5.00\\

chemequs & \hfill 0.67 & \hfill 0.61 & \hfill 0.57 & \hfill 0.6  &  \hfill 12.80 & \hfill 11.99\\

chemkin & \hfill 17.83 & \hfill 16.87 & \hfill 10.95 & \hfill 9.95 &  \hfill 32.34 & \hfill 29.83\\

cohn3 & \hfill 76.72 & \hfill 107.14 & \hfill 30.21 & \hfill 25.47 &  \hfill 37.73 & \hfill 39.20\\

cpdm5 & \hfill 1.78 & \hfill 1.57 & \hfill 1.69 & \hfill 1.68   & \hfill 0.69 & \hfill 0.70 \\

cyclic6 & \hfill 0.12 & \hfill 0.19 & \hfill 0.14 & \hfill 0.14  & \hfill 0.09 & \hfill 0.08 \\

cyclic7 & \hfill 58.72 & \hfill 60.94 & \hfill 68.59 & \hfill 65.28 & \hfill 6.64 & \hfill 7.08 \\

cyclic8 & \hfill 12056.24 & \hfill 14046.26 & \hfill 5826.18 & \hfill 4424.96 &  \hfill 235.73 & \hfill 245.65\\

d1 & \hfill 8.77 & \hfill 12.58 & \hfill 1.99 & \hfill 2.08 & \hfill 28.49 & \hfill 8.29\\

des18\_3 & \hfill 0.19 & \hfill 0.18 & \hfill 0.19 & \hfill 0.19 & \hfill 1.81 & \hfill 1.89\\

des22\_24 & \hfill 0.68 & \hfill 0.62 & \hfill 0.77 & \hfill 0.79 & \hfill 1.37 & \hfill 1.46\\

discret3 & \hfill 23322.8 & \hfill 20956.31 & \hfill 12642.49 & \hfill 13521.65 & \hfill 33658.09 & \hfill 19369.53\\

dl & \hfill 270.17 & \hfill 278.89 & \hfill 80.77 & \hfill 89.52 & \hfill 14.57 & \hfill 11.95\\

eco8 & \hfill 0.40 & \hfill 0.44 & \hfill 0.44 & \hfill 0.46 &  \hfill 0.20 & \hfill 0.20\\

eco9 & \hfill 3.22 & \hfill 5.60 & \hfill 4.99 & \hfill 5.08 & \hfill 1.25 & \hfill 1.20\\

eco10 & \hfill 52.56 & \hfill 56.70 & \hfill 65.71 & \hfill 68.06 &  \hfill 7.07 & \hfill 6.91\\

eco11 & \hfill 765.98 & \hfill 741.74 & \hfill 718.53 & \hfill 679.3 &  \hfill 62.33 & \hfill 51.08\\

%eco12 & \hfill 4083.00 & \hfill 8824.61 & \hfill 11224.11 & \hfill 11493.49 & \hfill (?) & \hfill 10416.81\\

extcyc5 & \hfill 1.35 & \hfill 1.53 & \hfill 1.46 & \hfill 1.37 &  \hfill 0.37 & \hfill 0.38\\

extcyc6 & \hfill 324.70 & \hfill 184.49 & \hfill 276.06 & \hfill 155.64 &  \hfill 45.36 & \hfill 47.96\\

extcyc7 & \hfill * & \hfill * & \hfill * & \hfill * &  \hfill 8242.00 & \hfill 8492.13\\

f744 & \hfill 4.88 & \hfill 7.71 & \hfill 2.22 & \hfill 2.68 & \hfill 1.47 & \hfill 1.38\\

f855 & \hfill 132.97 & \hfill 139.79 & \hfill 37.64 & \hfill 38.45 & \hfill 48.63 &  \hfill 37.06\\

fabrice24 & \hfill 108.52 & \hfill 116.77 & \hfill 8.2 & \hfill 7.7 & \hfill 9.45 & \hfill 8.70\\

filter9 & \hfill 20.97 & \hfill 5.76 & \hfill 1.13 & \hfill 1.6 & \hfill 80.04 &  \hfill 56.67\\

hairer2 & \hfill 62.91 & \hfill 108.17 & \hfill 126.69 & \hfill 125.43 & \hfill 92.07 &  \hfill 85.86\\

hairer3 & \hfill 1.96 & \hfill 0.92 & \hfill 0.32 & \hfill 1.4 & \hfill * &  \hfill *\\

hcyclic7 & \hfill 64.17 & \hfill 53.87 & \hfill 65.81 & \hfill 73.0 & \hfill 6.26 &  \hfill 6.76\\

hcyclic8 & \hfill 6024.97 & \hfill 4316.59 & \hfill * & \hfill 7560.99 & \hfill 229.70 &  \hfill 237.12\\

hf744 & \hfill 22.17 & \hfill 8.58 & \hfill 7.18 & \hfill 11.26 & \hfill 1.39 &  \hfill 1.32\\

hf855 & \hfill 2157.88 & \hfill 534.08 & \hfill 806.51 & \hfill 988.38 & \hfill 48.15 &  \hfill 36.69\\

hietarinta1 & \hfill 0.77 & \hfill 0.71 & \hfill 0.38 & \hfill 0.53 & \hfill 2.63 &  \hfill 2.15\\

i1 & \hfill 98.24 & \hfill 122.36 & \hfill 58.29 & \hfill 58.21 & \hfill 55.07 & \hfill 42.35\\

ilias13 & \hfill 1167.18 & \hfill 5851.97 & \hfill 3013.1 & \hfill 2469.62 & \hfill 336.21  & \hfill 309.64\\

ilias\_k\_2 & \hfill 323.59 & \hfill 669.68 & \hfill 445.51 & \hfill 270.21 & \hfill 55.41 &  \hfill 54.71\\

ilias\_k\_3 & \hfill 452.32 & \hfill 846.19 & \hfill 1162.7 & \hfill 622.14 & \hfill 90.67 &  \hfill 89.97\\

jcf26 & \hfill 224.96 & \hfill 211.24 & \hfill 16.44 & \hfill 14.65 & \hfill 31.64 & \hfill 25.59\\

katsura7 & \hfill 2.15 & \hfill 1.77 & \hfill 2.08 & \hfill 1.98 & \hfill 0.72 &  \hfill 0.79\\

katsura8 & \hfill 27.48 & \hfill 24.66 & \hfill 28.8 & \hfill 27.09 & \hfill 4.7 & \hfill 5.06\\

katsura9 & \hfill 337.52 & \hfill 294.59 & \hfill 340.45 & \hfill 311.98 & \hfill 33.47 &  \hfill 34.87\\

katsura10 & \hfill 4790.55 & \hfill 4983.11 & \hfill 7220.29 & \hfill 6204.95 & \hfill 287.38 & \hfill 292.02\\

kin1 & \hfill 15.18 & \hfill 20.32 & \hfill 7.11 & \hfill 7.11 & \hfill 50.56 & \hfill 45.33\\

kotsireas & \hfill 6.33 & \hfill 37.94 & \hfill 4.93 & \hfill 4.27 & \hfill 3.45 & \hfill 3.67\\

noon6 & \hfill 0.97 & \hfill 1.29 & \hfill 1.27 & \hfill 1.29 & \hfill 0.60 &  \hfill 0.62\\

noon7 & \hfill 28.87 & \hfill 32.58 & \hfill 37.52 & \hfill 38.52 & \hfill 4.93 & \hfill 4.77\\

noon8 & \hfill 1552.26 & \hfill 2292.84 & \hfill 3322.62 & \hfill 3152.57 & \hfill 43.65 & \hfill 42.80\\

pinchon1 & \hfill 10.37 & \hfill 0.04 & \hfill 0.01 & \hfill 0.01 & \hfill 4.09 & \hfill 3.54\\

rbpl & \hfill 210.94 & \hfill 177.51 & \hfill 173.8 & \hfill 173.98 & \hfill 38.33 & \hfill 35.79\\

rbpl24 & \hfill 108.78 & \hfill 116.78 & \hfill 8.23 & \hfill 7.7 & \hfill 9.62 & \hfill 8.74\\

redcyc6 & \hfill 0.16 & \hfill 0.17 & \hfill 0.13 & \hfill 0.14 & \hfill 0.10 & \hfill 0.10\\

redcyc7 & \hfill 913.75 & \hfill 1048.69 & \hfill 48.19 & \hfill 48.61 & \hfill 5.73 & \hfill 6.36\\

redeco10 & \hfill 18.51 & \hfill 18.66 & \hfill 23.91 & \hfill 22.4 & \hfill 2.33 & \hfill 2.40\\

redeco11 & \hfill 178.32 & \hfill 187.36 & \hfill 253.34 & \hfill 228.41 & \hfill 14.56 & \hfill 14.85\\

redeco12 & \hfill 1735.95 & \hfill 2172.75 & \hfill 4666.8 & \hfill 3385.97 & \hfill 101.51 & \hfill 103.02\\

reimer5 & \hfill 0.22 & \hfill 0.36 & \hfill 0.34 & \hfill 0.38 & \hfill 0.74 & \hfill 0.70\\

reimer6 & \hfill 9.69 & \hfill 21.60 & \hfill 24.19 & \hfill 23.96 & \hfill 42.13 &  \hfill 42.40\\

reimer7 & \hfill 719.37 & \hfill 3808.91 & \hfill 4756.4 & \hfill 4314.12 & \hfill 5216.53 & \hfill 5032.73\\

virasoro & \hfill 9.69 & \hfill 8.90 & \hfill 10.96 & \hfill 10.68 & \hfill 1.72 & \hfill  1.77
\end{longtable}
\end{scriptsize}
\end{center}

\section{Conclusion}

In this paper we experimentally investigated three different selection strategies for
the involutive algorithm specialized for Janet division and observed its computational
stability with respect to these variations in selection strategy. However, the problem of
finding heuristically best
selection strategies for the Janet division algorithms as well as for algorithms
exploiting other involutive division is still open and is of practical importance.
Our computer experimenting shows that the arbitrariness in selection of non-multiplicative
prolongation at the intermediate steps of the algorithm is sharply
grows with the number of variables and degree of the initial polynomials.

By this reason the number of prolongations with the same head degree or even
with the same leading monomial may achieve many hundreds and thousands for sufficiently
large examples. That is why, searching for heuristically best strategies is so important for
increasing computational efficiency of the involutive algorithmic methods.

\section{Acknowledgements}
The authors thank Daniel Robertz for running the benchmarks at RWTH, Aachen, with Magma V2.11-8 and V2.12-17.
The research presented in the paper was partially supported by
grants 04-01-00784 and 05-02-17645 from the Russian Foundation for Basic Research and by grant 2339.2003.2
from the Ministry of Education and Science of the Russian Federation.

\end{document}